\expandafter\chardef\csname pre amssym.def at\endcsname=\the\catcode`\@
\catcode`\@=11

\def\undefine#1{\let#1\undefined}
\def\newsymbol#1#2#3#4#5{\let\next@\relax
 \ifnum#2=\@ne\let\next@\msafam@\else
 \ifnum#2=\tw@\let\next@\msbfam@\fi\fi
 \mathchardef#1="#3\next@#4#5}
\def\mathhexbox@#1#2#3{\relax
 \ifmmode\mathpalette{}{\m@th\mathchar"#1#2#3}%
 \else\leavevmode\hbox{$\m@th\mathchar"#1#2#3$}\fi}
\def\hexnumber@#1{\ifcase#1 0\or 1\or 2\or 3\or 4\or 5\or 6\or 7\or 8\or
 9\or A\or B\or C\or D\or E\or F\fi}

\font\tenmsa=msam10
\font\sevenmsa=msam7
\font\fivemsa=msam5
\newfam\msafam
\textfont\msafam=\tenmsa
\scriptfont\msafam=\sevenmsa
\scriptscriptfont\msafam=\fivemsa
\edef\msafam@{\hexnumber@\msafam}
\mathchardef\dabar@"0\msafam@39
\def\dashrightarrow{\mathrel{\dabar@\dabar@\mathchar"0\msafam@4B}}
\def\dashleftarrow{\mathrel{\mathchar"0\msafam@4C\dabar@\dabar@}}

\def\ulcorner{\delimiter"4\msafam@70\msafam@70 }
\def\urcorner{\delimiter"5\msafam@71\msafam@71 }
\def\llcorner{\delimiter"4\msafam@78\msafam@78 }
\def\lrcorner{\delimiter"5\msafam@79\msafam@79 }
\def\yen{{\mathhexbox@\msafam@55 }}
\def\checkmark{{\mathhexbox@\msafam@58 }}
\def\circledR{{\mathhexbox@\msafam@72 }}
\def\maltese{{\mathhexbox@\msafam@7A }}

\font\tenmsb=msbm10
\font\sevenmsb=msbm7
\font\fivemsb=msbm5
\newfam\msbfam
\textfont\msbfam=\tenmsb
\scriptfont\msbfam=\sevenmsb
\scriptscriptfont\msbfam=\fivemsb
\edef\msbfam@{\hexnumber@\msbfam}

\catcode`\@=\csname pre amssym.def at\endcsname

\expandafter\ifx\csname pre amssym.tex at\endcsname\relax \else \endinput\fi
\expandafter\chardef\csname pre amssym.tex at\endcsname=\the\catcode`\@
\catcode`\@=11
\newsymbol\boxdot 1200
\newsymbol\boxplus 1201
\newsymbol\boxtimes 1202
\newsymbol\square 1003
\newsymbol\blacksquare 1004
\newsymbol\centerdot 1205
\newsymbol\lozenge 1006
\newsymbol\blacklozenge 1007
\newsymbol\circlearrowright 1308
\newsymbol\circlearrowleft 1309
\undefine\rightleftharpoons
\newsymbol\rightleftharpoons 130A
\newsymbol\leftrightharpoons 130B
\newsymbol\boxminus 120C
\newsymbol\Vdash 130D
\newsymbol\Vvdash 130E
\newsymbol\vDash 130F
\newsymbol\twoheadrightarrow 1310
\newsymbol\twoheadleftarrow 1311
\newsymbol\leftleftarrows 1312
\newsymbol\rightrightarrows 1313
\newsymbol\upuparrows 1314
\newsymbol\downdownarrows 1315
\newsymbol\upharpoonright 1316
 
\newsymbol\downharpoonright 1317
\newsymbol\upharpoonleft 1318
\newsymbol\downharpoonleft 1319
\newsymbol\rightarrowtail 131A
\newsymbol\leftarrowtail 131B
\newsymbol\leftrightarrows 131C
\newsymbol\rightleftarrows 131D
\newsymbol\Lsh 131E
\newsymbol\Rsh 131F
\newsymbol\rightsquigarrow 1320
\newsymbol\leftrightsquigarrow 1321
\newsymbol\looparrowleft 1322
\newsymbol\looparrowright 1323
\newsymbol\circeq 1324
\newsymbol\succsim 1325
\newsymbol\gtrsim 1326
\newsymbol\gtrapprox 1327
\newsymbol\multimap 1328
\newsymbol\therefore 1329
\newsymbol\because 132A
\newsymbol\doteqdot 132B
 
\newsymbol\triangleq 132C
\newsymbol\precsim 132D
\newsymbol\lesssim 132E
\newsymbol\lessapprox 132F
\newsymbol\eqslantless 1330
\newsymbol\eqslantgtr 1331
\newsymbol\curlyeqprec 1332
\newsymbol\curlyeqsucc 1333
\newsymbol\preccurlyeq 1334
\newsymbol\leqq 1335
\newsymbol\leqslant 1336
\newsymbol\lessgtr 1337
\newsymbol\backprime 1038
\newsymbol\risingdotseq 133A
\newsymbol\fallingdotseq 133B
\newsymbol\succcurlyeq 133C
\newsymbol\geqq 133D
\newsymbol\geqslant 133E
\newsymbol\gtrless 133F
\newsymbol\sqsubset 1340
\newsymbol\sqsupset 1341
\newsymbol\vartriangleright 1342
\newsymbol\vartriangleleft 1343
\newsymbol\trianglerighteq 1344
\newsymbol\trianglelefteq 1345
\newsymbol\bigstar 1046
\newsymbol\between 1347
\newsymbol\blacktriangledown 1048
\newsymbol\blacktriangleright 1349
\newsymbol\blacktriangleleft 134A
\newsymbol\vartriangle 134D
\newsymbol\blacktriangle 104E
\newsymbol\triangledown 104F
\newsymbol\eqcirc 1350
\newsymbol\lesseqgtr 1351
\newsymbol\gtreqless 1352
\newsymbol\lesseqqgtr 1353
\newsymbol\gtreqqless 1354
\newsymbol\Rrightarrow 1356
\newsymbol\Lleftarrow 1357
\newsymbol\veebar 1259
\newsymbol\barwedge 125A
\newsymbol\doublebarwedge 125B
\undefine\angle
\newsymbol\angle 105C
\newsymbol\measuredangle 105D
\newsymbol\sphericalangle 105E
\newsymbol\varpropto 135F
\newsymbol\smallsmile 1360
\newsymbol\smallfrown 1361
\newsymbol\Subset 1362
\newsymbol\Supset 1363
\newsymbol\Cup 1264
 
\newsymbol\Cap 1265
 
\newsymbol\curlywedge 1266
\newsymbol\curlyvee 1267
\newsymbol\leftthreetimes 1268
\newsymbol\rightthreetimes 1269
\newsymbol\subseteqq 136A
\newsymbol\supseteqq 136B
\newsymbol\bumpeq 136C
\newsymbol\Bumpeq 136D
\newsymbol\lll 136E
 
\newsymbol\ggg 136F
 
\newsymbol\circledS 1073
\newsymbol\pitchfork 1374
\newsymbol\dotplus 1275
\newsymbol\backsim 1376
\newsymbol\backsimeq 1377
\newsymbol\complement 107B
\newsymbol\intercal 127C
\newsymbol\circledcirc 127D
\newsymbol\circledast 127E
\newsymbol\circleddash 127F
\newsymbol\lvertneqq 2300
\newsymbol\gvertneqq 2301
\newsymbol\nleq 2302
\newsymbol\ngeq 2303
\newsymbol\nless 2304
\newsymbol\ngtr 2305
\newsymbol\nprec 2306
\newsymbol\nsucc 2307
\newsymbol\lneqq 2308
\newsymbol\gneqq 2309
\newsymbol\nleqslant 230A
\newsymbol\ngeqslant 230B
\newsymbol\lneq 230C
\newsymbol\gneq 230D
\newsymbol\npreceq 230E
\newsymbol\nsucceq 230F
\newsymbol\precnsim 2310
\newsymbol\succnsim 2311
\newsymbol\lnsim 2312
\newsymbol\gnsim 2313
\newsymbol\nleqq 2314
\newsymbol\ngeqq 2315
\newsymbol\precneqq 2316
\newsymbol\succneqq 2317
\newsymbol\precnapprox 2318
\newsymbol\succnapprox 2319
\newsymbol\lnapprox 231A
\newsymbol\gnapprox 231B
\newsymbol\nsim 231C
\newsymbol\ncong 231D
\newsymbol\diagup 231E
\newsymbol\diagdown 231F
\newsymbol\varsubsetneq 2320
\newsymbol\varsupsetneq 2321
\newsymbol\nsubseteqq 2322
\newsymbol\nsupseteqq 2323
\newsymbol\subsetneqq 2324
\newsymbol\supsetneqq 2325
\newsymbol\varsubsetneqq 2326
\newsymbol\varsupsetneqq 2327
\newsymbol\subsetneq 2328
\newsymbol\supsetneq 2329
\newsymbol\nsubseteq 232A
\newsymbol\nsupseteq 232B
\newsymbol\nparallel 232C
\newsymbol\nmid 232D
\newsymbol\nshortmid 232E
\newsymbol\nshortparallel 232F
\newsymbol\nvdash 2330
\newsymbol\nVdash 2331
\newsymbol\nvDash 2332
\newsymbol\nVDash 2333
\newsymbol\ntrianglerighteq 2334
\newsymbol\ntrianglelefteq 2335
\newsymbol\ntriangleleft 2336
\newsymbol\ntriangleright 2337
\newsymbol\nleftarrow 2338
\newsymbol\nrightarrow 2339
\newsymbol\nLeftarrow 233A
\newsymbol\nRightarrow 233B
\newsymbol\nLeftrightarrow 233C
\newsymbol\nleftrightarrow 233D
\newsymbol\divideontimes 223E
\newsymbol\varnothing 203F
\newsymbol\nexists 2040
\newsymbol\Finv 2060
\newsymbol\Game 2061
\newsymbol\mho 2066
\newsymbol\eth 2067
\newsymbol\eqsim 2368
\newsymbol\beth 2069
\newsymbol\gimel 206A
\newsymbol\daleth 206B
\newsymbol\lessdot 236C
\newsymbol\gtrdot 236D
\newsymbol\ltimes 226E
\newsymbol\rtimes 226F
\newsymbol\shortmid 2370
\newsymbol\shortparallel 2371
\newsymbol\smallsetminus 2272
\newsymbol\thicksim 2373
\newsymbol\thickapprox 2374
\newsymbol\approxeq 2375
\newsymbol\succapprox 2376
\newsymbol\precapprox 2377
\newsymbol\curvearrowleft 2378
\newsymbol\curvearrowright 2379
\newsymbol\digamma 207A
\newsymbol\varkappa 207B
\newsymbol\Bbbk 207C
\newsymbol\hslash 207D
\undefine\hbar
\newsymbol\hbar 207E
\newsymbol\backepsilon 237F
\catcode`\@=\csname pre amssym.tex at\endcsname

\magnification=1200
\hsize=468truept
\vsize=646truept
\voffset=-10pt
\parskip=4pt
\baselineskip=14truept
\count0=1

\dimen100=\hsize

\def\leftill#1#2#3#4{
\medskip
\line{$
\vcenter{
\hsize = #1truept \hrule\hbox{\vrule\hbox to  \hsize{\hss \vbox{\vskip#2truept
\hbox{{\copy100 \the\count105}: #3}\vskip2truept}\hss }
\vrule}\hrule}
\dimen110=\dimen100
\advance\dimen110 by -36truept
\advance\dimen110 by -#1truept
\hss \vcenter{\hsize = \dimen110
\medskip
\noindent { #4\par\medskip}}$}
\advance\count105 by 1
}
\def\rightill#1#2#3#4{
\medskip
\line{
\dimen110=\dimen100
\advance\dimen110 by -36truept
\advance\dimen110 by -#1truept
$\vcenter{\hsize = \dimen110
\medskip
\noindent { #4\par\medskip}}
\hss \vcenter{
\hsize = #1truept \hrule\hbox{\vrule\hbox to  \hsize{\hss \vbox{\vskip#2truept
\hbox{{\copy100 \the\count105}: #3}\vskip2truept}\hss }
\vrule}\hrule}
$}
\advance\count105 by 1
}
\def\midill#1#2#3{\medskip
\line{$\hss
\vcenter{
\hsize = #1truept \hrule\hbox{\vrule\hbox to  \hsize{\hss \vbox{\vskip#2truept
\hbox{{\copy100 \the\count105}: #3}\vskip2truept}\hss }
\vrule}\hrule}
\dimen110=\dimen100
\advance\dimen110 by -36truept
\advance\dimen110 by -#1truept
\hss $}
\advance\count105 by 1
}
\def\insectnum{\copy110\the\count120
\advance\count120 by 1
}

\font\ninerm=cmr9
\font\eightrm=cmr8

\font\tenrm=cmr10 at 10pt

\font\sc=cmcsc10

\def\msb{\fam\msbfam\tenmsb}

\def\bbc{{\msb C}}

\def\bbi{{\msb I}}

\def\bbq{{\msb Q}}
\def\bbr{{\msb R}}

\def\gra{\alpha}

\def\grd{\delta}
\def\gre{\epsilon}

\def\grl{\lambda}

\def\gro{\omega}

\def\grr{\rho}
\def\grs{\sigma}

\def\la#1{\hbox to #1pc{\leftarrowfill}}
\def\ra#1{\hbox to #1pc{\rightarrowfill}}

\def\fract#1#2{\raise4pt\hbox{$ #1 \atop #2 $}}

\def\bowtie{\hbox to 1pt{\hss}\raise.66pt\hbox{$\scriptstyle{>}$}
\kern-4.9pt\triangleleft}
\def\hsmash{\triangleright\kern-4.4pt\raise.66pt\hbox{$\scriptstyle{<}$}}
\def\boxit#1{\vbox{\hrule\hbox{\vrule\kern3pt
\vbox{\kern3pt#1\kern3pt}\kern3pt\vrule}\hrule}}

\def\za{\vrule height6pt width4pt depth1pt}

\font\aa=eufm10

\def\Got#1{\hbox{\aa#1}}

\def\cald{{\cal D}}

\def\calf{{\cal F}}

\def\calj{{\cal J}}

\def\call{{\cal L}}

\def\calz{{\cal Z}}

\def\ga{{\Got a}}

\def\gt{{\Got t}}

\def\gA{{\Got A}}

\def\gI{{\Got I}}

\def\gT{{\Got T}}

\font\svtnrm=cmr17

\font\bsc=cmcsc10 at 10truept

\def\kcon{2}
\def\acr{3}
\def\lcl{4}
\def\pr{5}
\def\cone{6}
\def\ee{7}

\def\Ric{\hbox{Ric}}

\def\cJ{\check{J}}
\def\cnabla{\check{\nabla}}
\def\endo{\hbox{End}}

\def\cs{\check{s}}

\centerline{\svtnrm Einstein Manifolds and Contact Geometry}
\bigskip
\centerline{\sc Charles P. Boyer~~ Krzysztof Galicki~~}
\footnote{}{\ninerm During the preparation of this work the authors 
were partially supported by NSF grant DMS-9970904.}
\bigskip
\centerline{\vbox{\hsize = 5.85truein
\baselineskip = 12.5truept
\eightrm
\noindent {\bsc Abstract}. We show that every K-contact Einstein
manifold is Sasakian-Einstein and discuss several corollaries of this
result.}}
\bigskip\bigskip\bigskip
\baselineskip = 10 truept
\centerline{\bf 1. Introduction} 
\bigskip

Recently the authors and their collaborators (cf. [BG1, BG2]) have used the
geometry of special types of Riemannian contact manifolds to construct
Einstein metrics of positive scalar curvature. In particular, in [BG2] we
studied Sasakian-Einstein geometry. Since Sasakian geometry is the odd
dimensional analogue of K\"ahler geometry, one might inquire as to the
validity of an odd dimensional Goldberg Conjecture. Recall that the
Goldberg conjecture [Gol] states that a compact almost K\"ahler manifold that
is also Einstein is K\"ahler-Einstein, that is the almost complex structure is
integrable. This conjecture has been confirmed by Sekigawa [Sek1, Sek2] in the
case of nonnegative scalar curvature. Since Sasakian-Einstein metrics
necessarily have positive scalar curvature, it is tempting to believe that an
odd dimensional Goldberg Conjecture hold true. The form in which one would
expect this conjecture to be true assumes that the metric be bundle-like. If
the Reeb vector field is quasi-regular so that under a compactness assumption
there is an orbifold fibration over an almost K\"ahler-Einstein orbifold it
seems quite likely that such a result should follow directly from Sekigawa's
result. However, in general one does not have such an orbifold submersion. We
handle this more general case by considering the closures of the leaves of the
characteristic foliation together with a construction of Molino [Mol1, Mol2]
which in the presence of a bundle-like Riemannian metric gives the existence
of a sheaf of commuting Killing vector fields. This sheaf can then be used to
approximate the geometry of the general case by orbifold submersions. Thus the
main purpose of this note is to prove the following:

\noindent{\sc Theorem} A: \tensl Let $(M,\eta,g)$ be a compact metric contact
manifold whose Riemannian metric $g$ is bundle-like with respect to the
characteristic foliation on $M.$ Then if $g$ is Einstein then it is
Sasakian-Einstein. Equivalently, every compact K-contact Einstein manifold is
Sasakian-Einstein. \tenrm

We also discuss some consequences of Theorem A to almost K\"ahler
structures on cones, and to some related work on $\eta$-Einstein manifolds.

\bigskip
\baselineskip = 10 truept
\centerline{\bf 2. Some Metric Contact Geometry} 
\bigskip

Let $(M,\cald)$ be a contact manifold and fix a contact 1-form $\eta$ such
that $\cald = \ker~\eta.$  The pair $(\cald,\gro)$, where $\gro$ is the
restriction of $d\eta$ to $\cald$ gives $\cald$ the structure of a symplectic
vector bundle. We denote by $\calj(\cald)$ the space of all almost complex
structures $J$ on $\cald$ that are compatible with $\gro,$ that is the subspace
of smooth sections $J$ of the endomorphism bundle $\endo~\cald$ that satisfy
$$ J^2= -\bbi, \qquad d\eta(JX,JY)=d\eta(X,Y), \qquad
d\eta(X,JX)>0\leqno{\kcon.1}$$ 
for any smooth sections $X,Y$ of $\cald.$ Notice
that each $J\in \calj(\cald)$ defines a Riemannian metric $g_\cald$ on $\cald$
by setting $g_\cald(X,Y) =d\eta(X,JY).$ One easily checks that $g_\cald$
satisfies the compatibility condition $g_\cald(JX,JY)=g_\cald(X,Y).$
Furthermore, the map $J\mapsto g_\cald$ is one-to-one, and the space
$\calj(\cald)$ is contractible. A choice of $J$ gives $M$ an almost CR
structure with a strictly pseudoconvex Levy form. 

Moreover, by extending $J$ to all of $TM$ one obtains an almost contact
structure [Bl,YK]. There are some choices of conventions to make here. We
define the section $\Phi$ of $\endo~TM$ by $\Phi =J$ on $\cald$ and
$\Phi\xi=0$, where $\xi$ is the Reeb vector field associated to $\eta.$ We can
also extend the transverse metric $g_\cald$ to a metric $g$ on all of $M$ by  
$$g(X,Y)= g_\cald +\eta(X)\otimes\eta(Y)= d\eta(X,\Phi Y)+
\eta(X)\otimes\eta(Y) \leqno{\kcon.2}$$ 
for all vector fields $X,Y$ on $M.$ One
easily sees that $g$ satisfies the compatibility condition $g(\Phi X,\Phi
Y)=g(X,Y)-\eta(X)\eta(Y).$ A contact manifold $M$ with a fixed contact form 
$\eta$ together with a vector field $\xi,$ a section $\Phi$ of $\endo~TM,$ and
a Riemannian metric $g$ which satisfy the conditions 
$$\eta(\xi)=1,\qquad \Phi^2=-\bbi +\xi\otimes \eta,\qquad g(\Phi X,\Phi Y)
=g(X,Y)-\eta(X)\eta(Y) \leqno{\kcon.3}$$
is known [Bl] as a {\it metric contact structure} on $M.$

Let us consider the characteristic foliation $\calf_\xi$ generated by the
Reeb vector field $\xi.$ If $\calf_\xi$ is a Riemannian foliation [Rei2, Mol1], 
that is, the holonomy pseudogroup induces isometries of Riemannian metrics on
the local leaf spaces of the local submersions defining $\calf_\xi.$ Then by
pulling back the metrics on the local leaf spaces one obtains a transverse
metric $g_\cald$ on the vector bundle $\cald$ that is invariant under the
(Reeb) flow generated by the Reeb vector field $\xi.$  This is equivalent to
the metric $g$ on $M$ given by \kcon.2 being bundle-like [Rei1, Rei2]. Hence

\noindent{\sc Definition} \kcon.4: \tensl A contact metric manifold
$(M,\eta,g)$ is said to be {\it bundle-like} if the Riemannian metric $g$ is
bundle-like. \tenrm

We have

\noindent{\sc Proposition} \kcon.5: \tensl On a complete contact metric
manifold $(M,\eta,g),$ the following are equivalent:
\item{(1)} $g$ is bundle-like.
\item{(2)} The Reeb flow is an isometry.
\item{(3)} The Reeb flow leaves the almost complex structure $J$ on $\cald$
invariant.
\item{(4)} The Reeb flow leaves the $(1,1)$ tensor field $\Phi$ invariant.
\item{(5)} The contact metric structure $(M,\eta,g)$ is K-contact.
\tenrm

\noindent{\sc Proof}: The conditions 
$$\eta(\xi)=1,\qquad \xi\rfloor d\eta =0$$
defining the Reeb vector field imply that both the symplectic form $d\eta$ and
the contact form $\eta$ are invariant under the Reeb flow.  From its definition
[Rei1] $g$ is bundle-like if and only if the transverse metric $g_\cald$ is
basic, that is, if and only if the Reeb flow leaves $g_\cald$ invariant as
well. Since $d\eta$ is invariant under the Reeb flow, $g_\cald$ is invariant
if and only if $J$ is invariant, or equivalently if and only if $\Phi$ is
invariant. \hfill\za

\noindent{\sc Remarks} \kcon.6: We prefer the appellation bundle-like contact
metric structure to the more common K-contact structure, since it is more
descriptive and emphasizes the foliation aspect. We shall use these two
terms interchangeably depending on the context. We also refer to the transverse
structure $(d\eta,J,g_\cald)$ on $M$ as a {\it transverse almost K\"ahler
structure}.

There are obstructions to admitting K-contact metric structures. Indeed, it is
well known [Bl] that K-contact metrics on a $2n+1$ dimensional manifold can be
characterized by the condition that the Ricci tensor equals $2n$ in the
direction of the Reeb vector field $\xi.$ Thus, any metric of non-positive
Ricci curvature cannot have a K-contact metric in its homothety class. 
However, it is much stronger to obtain obstructions which only depend on the
smooth structure of the manifold.  We mention one such result that follows
directly from the work of Gromov [Gr], Carriere [Car], and Inoue and Yano
[IY]. 

\noindent{\sc Theorem} \kcon.7: \tensl If a compact manifold $M$ admits a
bundle-like contact metric structure, then the Gromov invariant $||M||$ and
all the Pontrjagin numbers of $M$ vanish. In particular, if a compact manifold
$M$ admits a decomposition as a connected sum $M=M_1\#\cdots \# M_k$, where for
some $i=1,\cdots,k$ the manifold $M_i$ admits a metric of strictly negative
sectional curvature, then $M$ does not admit any bundle-like contact metric
structure. \tenrm 

\bigskip
\baselineskip = 10 truept
\centerline{\bf 3. The Almost CR Structure} 
\bigskip

In this section we consider the integrability of the almost CR structure $J.$
Let $(M,\eta,g)$ be a contact metric manifold. The almost CR structure $J$ is
integrable, that is, $(\cald,J)$ defines a CR-structure on $M$ if and only if
for any smooth sections $X,Y$  of $\cald$ the following conditions hold:
\item{(i)} $[X,JY]+[JX,Y]$ is a smooth section of $\cald.$ 
\item{(ii)} $J[X,JY]+J[JX,Y]= [JX,JY]-[X,Y].$ 

\noindent In our case condition (i) follows automatically from the
antisymmetry of the symplectic form $d\eta.$ Condition (ii) is the vanishing of
the Nijenhuis tensor of $J.$ 

Now let $\nabla$ denote the Levi-Civita connection with respect to the metric
$g$ on $M.$ By restricting to $\cald$ and taking the horizontal projection we
get an induced connection  $\nabla^\cald$ on $\cald$ defined by [Ton]
$$\nabla^\cald_XY = \cases{(\nabla_XY)^h&if $X$ is a smooth section of $\cald$,
\cr                         [\xi,Y]^h&if $X=\xi$, \cr}    \leqno{\acr.1}$$
where $Y$ is a smooth section of $\cald$ and the superscript $h$ denotes the
projection onto $\cald.$ An entirely standard computation gives

\noindent{\sc Proposition} \acr.2: \tensl Let $(M,\eta,g)$ be an contact metric
manifold. Then $\nabla^\cald J=0$ if and only if the almost CR-structure $J$ on
$\cald$ is integrable and $\call_\xi\Phi=0.$ \tenrm

\noindent{\sc Proof}: First we notice that as mentioned above 
condition (ii) above is automatically satisfied. Next, one easily sees that
the invariance of $\Phi$ under $\xi$ holds if and only if $\nabla^\cald_\xi
J=0.$ Now the connection $\nabla^\cald$ is torsion-free [Ton], so 
$$\nabla_X^\cald Y-\nabla_Y^\cald X=[X,Y]^h.$$
Now as above the vertical part of $N_J$ vanishes, and a straightforward
computation gives 
$$N_J(X,Y)^h= (\nabla^\cald_{JX}J)(Y)-(\nabla^\cald_{JY}J)(X)
+J(\nabla^\cald_YJ)(X)  -J(\nabla^\cald_XJ)(Y).$$
The only if part clearly holds (which is all we shall need), and the if
part follows by a standard computation (c.f. [YK]). \hfill\za

Recall that a contact metric structure $(M,\eta,g)$ is said to be {\it normal}
if the Nijenhuis tensor $N_\Phi$ defined by
$$N_{\Phi}(X,Y) ~=~ [\Phi X,\Phi Y]+(\Phi)^2 [X,Y]- \Phi [X,\Phi Y] 
 - \Phi [\Phi X,Y] \leqno{\acr.3}$$ 
satisfies
$$N_{\Phi} = -d\eta\otimes \xi. \leqno{\acr.4}$$ 
A normal contact metric structure on $M$ is also called a {\it Sasakian}
structure.

\noindent{\sc Proposition} \acr.5: \tensl Let $(M,\eta,g)$ be an contact metric
manifold. Then $(M,\eta,g)$ is normal (Sasakian) if and only if the almost
CR-structure $J$ is integrable and $\call_\xi\Phi=0.$  \tenrm 

\noindent{\sc Proof}: For any vector fields $X,Y$ on $M$ we have 
$$N_{\Phi}(X,Y)+d\eta(X,Y)\xi = [\Phi X,\Phi Y]+(\Phi)^2 [X,Y]- \Phi 
[X,\Phi Y] - \Phi [\Phi X,Y] +d\eta(X,Y)\xi. \leqno{\acr.6}$$ 
If $X$ and $Y$ are both horizontal then this equals 
$$[\Phi X,\Phi Y]-[X,Y] -\Phi [\Phi X,Y]-\Phi [X,\Phi Y]$$   
whose vanishing is equivalent to (iii) above. Also applying $\eta$ to this
equation and replacing $X$ by $\Phi X$ implies (ii). If one  vector field is
vertical, say $X=\xi$ then we have 
$$N_{\Phi}(\xi,Y)+d\eta(\xi,Y)\xi = (\xi\rfloor d\eta)(Y)\xi -\Phi\circ 
(\call_{\xi}\Phi)(Y).$$ 
So the result follows. \hfill\za 

\bigskip
\centerline{\bf 4. The Leaf Closures of $\calf_\xi$} 
\bigskip

In this section we study the leaf closures of the characteristic foliation. In
[Mol1, Mol2] Molino has shown that on any compact connected manifold $M$ with a
Riemannian foliation $\calf$ there is a locally constant sheaf $C(M,\calf)$,
called the commuting sheaf, consisting of germs of local transverse vector
fields that are Killing vector fields with respect to the transverse metric,
and whose orbits are precisely the closures of the leaves of $\calf.$ 
Moreover, Carri\`ere [Car] (See also the appendix in [Mol1] by Carri\`ere) has
shown in the case of Riemannian foliations of dimension one (Riemannian flows)
that the leaf closures are diffeomorphic to tori, and that the flow is
conjugate by the diffeomorphism to a linear flow on the torus.  

Here we adapt this to our situation, that is, $(M,\eta,g)$ is a compact
bundle-like metric contact manifold. We denote the isometry group of $(M,g)$
by $\gI(M,g),$ and the group of automorphisms of the K-contact structure
$(M,\eta,g),$ by $\gA(M,\eta,g).$  When $M$ is compact the well known
theorem of Myers and Steenrod says that $\gI(M,g)$ is a compact Lie group.
Moreover, $\gA(M,\eta,g)$ is a closed Lie subgroup of $\gI(M,g)$ [Tan1].  In
our case the Reeb flow belongs to the automorphism group $\gA(M,\eta,g)$ which
is a compact Lie group. Thus, the closure $\gT$ of the Reeb flow is a compact
commutative Lie group, i.e., a torus, which lies in $\gA(M,\eta,g).$ Now the
Reeb flow is a strict contact transformation lying in the center of the group
of strict contact transformations [LM]; hence, it lies in the center of 
$\gA(M,\eta,g).$ It follows that $\gT$ also lies in the center of
$\gA(M,\eta,g).$  Summarizing  we have

\noindent{\sc Proposition} \lcl.1: \tensl Let $(M,\eta,g)$ be a compact
bundle-like contact metric manifold. Then the leaf closures of the Reeb flow
are the orbits of a torus $\gT$ lying in the center of
the Lie group $\gA(M,\eta,g)$ of automorphisms of $(M,\eta,g),$ and the Reeb
flow is the orbit of a linear flow on $\gT.$   \tenrm   

The dimension of the torus in Proposition \lcl.1 is an invariant of the
K-contact structure that we call the {\it rank} of $(M,\eta,g)$ and denote
by $\hbox{rk}(M,\eta).$  We have (see also [Ruk2])

\noindent{\sc Lemma} \lcl.2: \tensl Let $(M,\eta,g)$ be a compact
bundle-like contact metric manifold of dimension $2n+1.$ Then the rank
$\hbox{rk}(M,\eta)$ depends only on the Pfaffian structure $(M,\eta)$ and 
satisfies $1\leq\hbox{rk}(M,\eta)\leq n+1.$ \tenrm

\noindent{\sc Proof}: Consider the Lie algebra $\gt$ of $\gT.$ It consists of
the Reeb vector field $\xi$ together with the infinitesimal generators of the
leaf closures. The projections of these generators onto $\cald$ are global
sections of Molino's commuting sheaf $C(M,\calf_\xi).$ Thus, they give integral
submanifolds of the subbundle $\cald.$ It is well known [LM] that the integral
submanifolds of maximal dimension, that is the Legendre submanifolds of the
contact structure, have dimension $n.$ Hence, these together with the Reeb
vector field generate a torus of dimension at most $n+1.$ Furthermore, Molino
[Mol1] shows that the commuting sheaf is independent of the transverse metric,
so $\hbox{rk}(M,\eta)$ is independent of $g.$ \hfill\za

Now the rank $\hbox{rk}(M,\eta)$ is not an invariant of the contact
structure $(M,\cald)$ but only of the Pfaffian structure $(M,\eta).$ The case
$\hbox{rk}(M,\eta)=1$ is the quasi-regular case, while the other extreme 
$\hbox{rk}(M,\eta)=n+1$ is the toric case studied in [BM1, BM2, BG3].
Furthermore, Rukimbira [Ruk1] showed that one can approximate any K-contact
form $\eta$ by a sequence of quasi-regular K-contact forms in the same contact
structure. Thus, every K-contact manifold has an $\eta$ of rank 1. Since we
shall discuss this approximation in detail in the next section, we only mention
here that one chooses a sequence of vector fields $\xi_j$ in $\gt$ with
periodic orbits that converges to the Reeb vector field $\xi.$ Then the dual
1-forms $\eta_j$ are quasi-regular contact forms in the same contact
structure.

\bigskip
\centerline{\bf 5. The Proof of Theorem A} 
\bigskip

We shall prove the following restatement of Theorem A:

\noindent{\sc Theorem} A$'$: \tensl Let $(M,\eta,g)$ be a  compact
K-contact Einstein manifold. Then $(M,\eta,g)$ is Sasakian-Einstein. \tenrm

\noindent{\sc Proof}: We first prove the theorem under the assumption that 
$\eta$ is quasi-regular. By Thomas [Tho] and [BG1], $M$ is the total space of
a principal $S^1$ V-bundle over a compact almost K\"ahler orbifold  $\calz.$
Furthermore, by [Bes] the induced metric $h$ on $\calz$ is almost
K\"ahler-Einstein which has positive scalar curvature, since $g$ has  positive
scalar curvature. Now since Sekigawa's [Sek1, Sek2] proof of the  Goldberg
conjecture in the case of nonnegative scalar curvature only involves local
curvature computations together with a Bochner type argument using Stokes
Theorem, it carries over to the case of a compact orbifold.  So the almost
complex structure on $\calz$ is integrable, and $(\calz,h)$ is
K\"ahler-Einstein. It then follows from the orbifold version of Hatakeyama
[Hat] that $(M,\eta,g)$ is normal, hence, Sasakian-Einstein. This proves
the result under the assumption of quasi-regularity. 

Now assume that $(M,\eta,g)$ is K-contact and Einstein, but not quasi-regular.
Then by Proposition \lcl.1 the Reeb vector field $\xi$ lies in the commutative
Lie subalgebra $\gt(M,\calf_\xi)\subset \ga(M,g)$ which has dimension
$k>1.$ Thus, there exists a sequence of quasi-regular contact forms $\eta_j$
and Reeb vector fields $\xi_j\in \gt(M,\calf_\xi)$ that approximate
$(\eta,\xi)$ in the compact-open $C^{\infty}$ topology. (In what follows we use
this topology on the space of smooth sections of all tensor bundles.)
Explicitly, there is a monotonically decreasing sequence $\{\gre_j\}_1^\infty$
with ${\displaystyle{\lim_{j\to\infty}}}\gre_j=0$ such that                 
$$\eta_j =f(\gre_j)\eta, \qquad \xi_j= \xi+\grr_j, \qquad f(\gre_j)=  {1\over
1+\eta(\grr_j)},\leqno{\pr.1}$$ 
where  $f(\gre_j)$ are positive functions in
$C^\infty(M)$ that  satisfy ${\displaystyle{\lim_{j\to\infty}}}f(\gre_j)=1.$ 
Clearly $\grr_j\in \gt(M,\calf_\xi)$ and
${\displaystyle{\lim_{j\to\infty}}}\grr_j=0.$ Moreover, $\hbox{ker}~\eta_j
=\hbox{ker}~\eta=\cald,$ so we have the same underlying contact structure.  We
also have the following easily verified relations for the induced contact
endomorphisms $\Phi_j:$ $$\Phi_j=\Phi -{1\over
1+\eta(\grr_j)}\Phi\grr_j\otimes \eta  =\Phi -f(\gre_j)\Phi\grr_j\otimes
\eta.\leqno{\pr.2}$$ This implies that $\Phi_j\xi_j=0$ and that the almost
complex structure $J$ on $\cald$ remains unchanged. However, the induced
metrics become $$g_j= f(\gre_j)g_\cald \oplus f(\gre_j)^2\eta\otimes \eta =g-
\eta(\grr_j)\bigl(g_\cald +2\eta\otimes \eta\bigr)+o(\gre_j^2).$$ For $\gre_j$
small enough $g_j$ are well defined Riemannian metrics on $M$ which can easily
be seen to satisfy the compatibility conditions
$$g_j(\Phi_jX,\Phi_jY)=g_j(X,Y)-\eta_j(X)\eta_j(Y).\leqno{\pr.3}$$  Moreover,
since $\xi_j\in \gt\subset \ga(M,\eta),$ it follows that the functions
$f(\gre_j)\in C^\infty(M)^\gT$, where $C^\infty(M)^\gT$ denotes the subalgebra
of $C^\infty(M)$ invariant under the action of the torus $\gT.$ Thus, from
\pr.2 we have  $$\call_{\xi_j}\Phi_j=0.\leqno{\pr.4}$$
Hence $(M,\eta_j,\xi_j,\Phi_j,g_j)$
is a sequence of quasi-regular K-contact structures on $M$ whose limit with
respect to the compact-open $C^\infty$ topology is the original K-contact
Einstein structure $(M,\eta,\xi,\Phi,g).$ Now the metrics $g_j$ are not
Einstein, but their Ricci tensor can be seen to satisfy 
$$\hbox{Ric}_{g_j}=\grl_jg_j +A(\gre_j,\grr_j,g),\leqno{\pr.5}$$
where $A(\gre_j,\grr_j,g)$ is a traceless symmetric 2-tensor field depending on
$\gre_j,\grr_j,g$ that tends to $0$ with $\gre_j,$ and $\grl_j\in
C^\infty(M)$ satisfy ${\displaystyle{\lim_{j\to\infty}}}\grl_j =2n.$ 

Now there is a sequence of orbifold Riemannian submersions $\pi_j:M\ra{1.5}
\calz_j$, where $(\calz_j,h_j)$ are a sequence of compact almost K\"ahler
orbifolds satisfying $\pi^*_jh_j= f(\gre_j)g_\cald.$ Moreover, it follows from
the above limits that the scalar curvatures of the $h_j$ are all positive.
Notice that in Sekigawa's proof [Sek2] of the positive scalar curvature
Goldberg conjecture, the Einstein condition is not used until section 4 of
[Sek2]. Following [Sek2] and making the necessary adjustments to our
situation, we find that there are nonnegative numbers $\grd_j$ and nonnegative
smooth functions $F_j$ such that 
$$\int_{\calz_j}\bigl(F_j+{\cs_j\over n}||\cnabla_j\cJ_j||_{\calz_j}^2+
{1\over 2n}||\cnabla_j\cJ_j||_{\calz_j}^4\bigr)\grs_j\leq \grd_j,\leqno{\pr.6}$$
where $\cnabla_j,\cJ_j,\cs_j,\grs_j$ and $||\cdot||_{\calz_j}$ are the
Levi-Civita connection,  almost complex structure, scalar curvature, volume
element, and Riemannian norm, respectively on $(\calz_j,h_j).$ Now since the
metrics $g,g_j$ are bundle-like the leaves of the characteristic foliation are
geodesics and the O'Neil tensors $T$ and $N$ vanish [Ton]. Moreover, for
any K-contact manifold of dimension $2n+1$ the O'Neill tensor $A$ satisfies
$||A||^2=g(A\xi,A\xi)=2n.$ Thus, we have the relation between the functions
$\grl_j$ on $M$ and the scalar curvatures $\cs_j$ on $\calz_j:$         
$$\cs_j =(2n+1)\grl_j +2n. \leqno{\pr.7}$$                                   
So that ${\displaystyle{\lim_{j\to\infty}}}\cs_j=
2n+(2n+1){\displaystyle{\lim_{j\to\infty}}}\grl_j=
4n(n+1).$ Furthermore, we have ${\displaystyle{\lim_{j\to\infty}}}\grd_j=0.$ 
Thus, since
$F_j$ (see [Sek2]) and $\cs_j$ are nonnegative for each $j,$ the estimate
\pr.6 implies the estimate  $$||\cnabla_j\cJ_j||_{\calz_j}\leq
\grd'_j,\leqno{\pr.8}$$ where $\grd'_j$ are nonnegative numbers satisfying 
${\displaystyle{\lim_{j\to\infty}}}\grd'_j=0.$ Now for each $j$ the horizontal
lift of $\cnabla\cJ_j$ is the horizontal projection $(\nabla_jJ_j)^h= (\nabla_j
\Phi_j)^h$, where $\nabla_j,J_j,\Phi_j$ are the corresponding Levi-Civita
connection and tensor fields with respect to the metrics on $M.$ But on $M$
$J_j=J$ for all $j$ and we have
$$||(\nabla J)^h||= {\displaystyle{\lim_{j\to\infty}}}||(\nabla_jJ)^h||_j \leq
{\displaystyle{\lim_{j\to\infty}}}\grd'_j =0,\leqno{\pr.9}$$ 
where $||\cdot ||_j$ is the Riemannian norm with respect to $g_j.$ So by 
Proposition \acr.2 the almost CR structure on $\cald$ is integrable which by
Proposition \acr.5 implies that $(M,\eta,g)$ is Sasakian-Einstein.
\hfill\za 

\noindent{\sc Remark} \pr.10: Notice that by \pr.7 the scalar
curvatures $\cs_j$ of the orbifolds $(\calz_j,h_j)$ are close to $4n(n+1).$
But the integrability argument actually holds for a much larger range of
scalar curvatures, namely $\cs_j\geq 0.$ We shall make use of this in section
\ee~ when discussing $\eta$-Einstein metrics.   

\bigskip
\centerline{\bf 6. Almost K\"ahler Cones} 
\bigskip

Here we give a corollary of Theorem A concerning almost K\"ahler cones. 
We consider the symplectification of $(M,\eta),$ namely, the symplectic
cone $$C(M)=(M\times \bbr^+,d(r^2\eta)).$$ We can extend the almost complex
structure $J$ on $\cald$ to an almost complex structure $I$ on $TC(M)$ by 
setting 
$$ I=J ~\hbox{on}~ \cald, \qquad I\xi=-\Psi,\qquad I\Psi =\xi,\leqno{\cone.1}$$
where $\Psi= r{\partial \over \partial r}$ is the Euler vector field. Then
the metric $g_\cald+\eta\otimes \eta$ on $M$ corresponds to the metric
$dr^2+r^2(g_\cald+\eta\otimes \eta)$ on $C(M).$ We have arrived at

\noindent{\sc Proposition} \cone.2: \tensl Contact metric geometry on $M$ 
corresponds to almost K\"ahler geometry on $C(M).$ \tenrm

So what does the K-contact condition on $M$ correspond to? We have

\noindent{\sc Proposition} \cone.3: \tensl A compact metric contact manifold
$(M,\eta,g_M)$ is K-contact if and only if $(C(M),d(r^2\eta),dr^2+r^2g_M)$ is
almost K\"ahler with $\Psi-i\xi$  pseudo-holomorphic. \tenrm

\noindent{\sc Proof}: By the properties of the Reeb vector field one
easily sees that $\Psi-i\xi$ is pseudo-holomorphic, i.e., an infinitesimal
automorphism of $I$ if and only if $\call_\xi J=0.$ But this holds if and 
only if $\call_\xi g_\cald=0.$
\hfill\za
 
\noindent{\sc Remark} \cone.4: In the case that the complex vector field
$\Psi-i\xi$ on $C(M)$ is not pseudo-holomorphic, a quotient formed
by dividing by the resulting $\bbc^*$ action, or equivalently, by the
symplectic reduction of the $S^1$ action will lose both the almost complex
structure and the Riemannian structure.

\noindent{\sc Corollary} \cone.5: \tensl Let $M$ be compact with a metric
contact structure $(\eta,g_M),$ and consider the almost K\"ahler cone
$(C(M),d(r^2\eta),dr^2+r^2g_M).$ Suppose that the $(1,0)$ vector field
$\Psi-i\xi$ is pseudo-holomorphic and the cone metric $dr^2+r^2g_M$ is Ricci
flat, then the almost complex structure is integrable and the cone metric is
Calabi-Yau. \tenrm 

\noindent{\sc Proof}: It is well known that a cone metric $dr^2+r^2g_M$ on a
cone $C(M)$ of dimension $N$ is Ricci flat if and only if the metric $g_M$ is
Einstein with Einstein constant $N-2.$ Thus, the result follows from
Proposition \cone.3 and Theorem A. \hfill\za

Another easy consequence of our results involve Vaisman's generalized Hopf
manifolds. Consider an almost K\"ahler cone $(C(M),dr^2+r^2g_M).$ Then on
$M\times S^1$ defined as the quotient manifold of $C(M)$ by the discrete group
generated by $r\mapsto e^ar$, where $0<a<1$ is fixed, the metric $g_M+({dr\over
r})^2$ is locally conformally almost K\"ahler. Furthermore, the vector field
$\Psi$ and the almost complex structure $I$ pass to the quotient. Then we have

\noindent{\sc Corollary} \cone.6: \tensl Let $M$ be compact with a metric
contact structure $(\eta,g_M)$ and consider the locally conformally almost
K\"ahler manifold $(M\times S^1,g_M+({dr\over r})^2).$ Suppose further that the
$(1,0)$ vector field $\Psi -i\xi$ is pseudo-holomorphic and that the locally
defined almost K\"ahler metrics $dr^2+r^2g_M$ are Ricci flat. Then the almost
complex structure $I$ is integrable so the manifold  $(M\times
S^1,g_M+({dr\over r})^2)$ is a locally conformal Calabi-Yau manifold [BG2].
\tenrm

\bigskip
\centerline{\bf 7. Some Remarks on $\eta$-Einstein Metrics} 
\bigskip
 
We conclude with some results about $\eta$-Einstein metrics. First
recall [Tan2,YK]

\noindent{\sc Definition} \ee.1: \tensl A metric contact structure $(\eta,g)$
on $M$ is said to be $\eta$-Einstein if there are constants $a,b$ such that
$\Ric_g= ag+b\eta\otimes \eta.$ \tenrm

Actually if $(M,\eta,g)$ is Sasakian and such a condition
holds for $\Ric_g$, where $a,b$ are smooth functions, then these functions must
be constant [YK]. In this section we shall prove:

\noindent{\sc Theorem} \ee.2: \tensl Let $(M,\eta,g)$ be a compact K-contact
manifold such that $g$ is $\eta$-Einstein. Then
\item{(i)} If $a>-2$ the almost CR-structure $J$ is integrable, so $g$ is
Sasakian. Moreover, for $\gra= \displaystyle{{a+2\over 2n+2}}$ the metric $\gra
g+ \gra(\gra-1)\eta\otimes \eta$ is Sasakian-Einstein. Hence, $\pi_1(M)$ is
finite.
\item{(ii)} If $a=-2$ the almost CR-structure is integrable, so $(g,\eta)$ is
Sasakian $\eta$-Einstein. Moreover, if $M$ has finite fundamental group then
$\hbox{rk}(M,\eta)=1$ so the K-contact manifold $(M,\eta,g)$ is
quasi-regular, and the total space of a principal $S^1$ V-bundle over a
Calabi-Yau orbifold.
\item{(iii)} If $a<-2$ then $\hbox{rk}(M,\eta)=1$ so the K-contact manifold is
quasi-regular, and the total space of a principal $S^1$ V-bundle over an
almost K\"ahler-Einstein orbifold with Einstein constant $2n(a+2).$  \tenrm

\noindent{\sc Proof}: (i): Notice that Tanno [Tan2] proves the second statement
of (i) under the assumption that $(M,\eta,g)$ is Sasakian. However, as we shall
see this assumption is not necessary. Since the O'Neill tensors $T$ and $N$
vanish and $A$ satisfies $g(A_X,A_Y)=g(\Phi X,\Phi Y)=g(X,Y),$ it follows that
the Ricci curvature of the transverse metric $g_\cald$ satisfies
$$\Ric_{g_\cald}=\Ric_g|_{\cald\times \cald} + 2g|_{\cald\times \cald}. 
\leqno{\ee.3}$$
The condition that on $\cald$ the Ricci curvature satisfies $\Ric_g >-2$ is
equivalent to the condition that $\Ric_{g_\cald} >0.$ Now even though in
general we do not have a Riemannian submersion (even in the orbifold sense),
the canonical variation described in Besse [Bes] applies equally well to our
foliation since it is based on the O'Neill formulas which do hold in our case.
Then one easily sees that by choosing $\gra= \displaystyle{{a+2\over 2n+2}}$
the metric $g' =\gra g+ \gra(\gra-1)\eta\otimes \eta$ is Einstein, so
$(M,\gra\eta,g')$ is K-contact and Einstein. Thus, by Theorem A it is
Sasakian-Einstein. Since the underlying almost CR-structure hasn't changed,
the original K-contact structure $(M,\eta,g)$ is Sasakian; hence, it is
Sasakian $\eta$-Einstein.

\noindent (iii): Let $a<-2.$ Then from \ee.3 $\Ric_{g_\cald} <0.$ Suppose that
the K-contact structure $(M,\eta,g)$ is not quasi-regular. Then Molino's
commuting sheaf $C(M,\calf_\xi)$ is non-vanishing. So by a perturbation of the
K-contact structure there is a quasi-regular K-contact structure $(\eta',g')$
with the same commuting sheaf. But by a theorem of Molino and Sergiescu [MoSe]
the sheaf $C(M,\calf_\xi)$ has a global trivialization. Thus, there are
transverse Killing vector fields on $(M,\eta',g').$ These project to
non-trivial Killing vector fields on a compact orbifold $\calz'$ with negative
Ricci curvature. But as in the manifold case a compact orbifold with negative
Ricci curvature can have no Killing fields. This gives a contradiction.

\noindent (ii): By \ee.3 the case $a=-2$ corresponds to the vanishing of the
scalar curvature of the transverse metric $g_\cald.$ But then as mentioned in
Remark \pr.10 the proof of Theorem A holds in this case, and the almost
CR-structure $J$ is integrable. In this case the transverse geometry is
Calabi-Yau, and the metric $g$ is Sasakian $\eta$-Einstein with $a=-2$ and
$b=2n+2.$ To prove the second statement we proceed as in the proof of (iii),
only now the transverse Ricci tensor vanishes implying that any Killing fields
on $\calz'$ must be parallel. Moreover, the vector space of these Killing
fields has dimension equal to the first Betti number $b_1(\calz').$ But the
finiteness of $\pi_1(M)$ together with the long exact homotopy sequence of the
orbifold fibration $\pi:M\ra{1.3} \calz'$ implies that $\pi_1^{orb}(\calz')$
is also finite, and this implies that
$H_1(\calz',\bbq)=H_1^{orb}(\calz',\bbq)=0.$ This gives a contradiction. 
\hfill\za

Actually our proof gives a bit more:

\noindent{\sc Proposition} \ee.4: \tensl Let $(M,\eta,g)$ be a compact
K-contact manifold. Suppose that the transverse Ricci tensor satisfies
$\Ric_{g_\cald}\leq 0$ and $\Ric_{g_\cald}(v,v)<0$ for all $v\in
\cald_p-\{0\}$ for some $p\in M.$ Then $\hbox{rk}(M,\eta)=1$ so the
K-contact manifold is quasi-regular, and the total space of a principal $S^1$
V-bundle over an almost K\"ahler orbifold. \tenrm

\bigskip
\centerline{\bf 8. A Remark on Contact 3-Structures} 
\bigskip

It is interesting to inquire about a quaternionic analogue of our main theorem.
In this regard it has been recently observed by Kashiwada [Kas] that a much
stronger result is available. Indeed, a {\it metric contact 3-structure} on
a manifold $M$ is a triple of contact structures
$\{\eta^a,\xi^a,\Phi^a\}_{a=1}^3$ associated with the same metric $g$ such
that
$$\Phi^a\circ \Phi^b-\xi^a\otimes \eta^b ~=~  -\gre^{abc}\Phi^c
-\grd^{ab}\hbox{id}$$
where $\gre^{abc}$ is the totally antisymmetric symbol and sum over repeated
indices is used. If each contact structure is normal the triple is called {\it
3-Sasakian} (cf. [BG1]). 

\noindent{\sc Theorem} [Kashiwada]: \tensl Every metric contact 3-structure is
3-Sasakian. \tenrm

The key to the proof of this theorem is a result of Hitchin buried deep in his
famous stable pairs paper [Hit]. This result says that an almost hyperk\"ahler
structure must be hyperk\"ahler. More explicitly, if one has a manifold with a
triple of almost complex structures satisfying the algebra of the quaternions,
together with a triple of compatible K\"ahler forms all of which are closed,
then the almost complex structures are integrable. That is, the quaternionic
algebra and the closedness of the forms are strong enough to force
integrability. Then Kashiwada's Theorem follows from Hitchin's Lemma together
with the following quaternionic analogue of Proposition \cone.2:

\noindent{\sc Proposition} 8.1: \tensl A Riemannian manifold $(M,g)$ has a
compatible contact 3-structure $(\eta^a,\xi^a,\Phi^a)$ if and only if the cone
$(C(M),dr^2+r^2g)$ is almost hyperk\"ahler. Furthermore, $(M,g)$ is 3-Sasakian
if and only if $(C(M),dr^2+r^2g)$ is hyperk\"ahler. \tenrm

We conclude by mentioning that some 
weaker results have been obtained in the last
few years, first by Tanno [Tan3] and then by Jelonek [Jel]. 
In 1996 Tanno observed that
in dimension 7 any 3-K-contact manifold must be 3-Sasakian and later
Jelonek (using similar
techniques but assuming leaf compactness of the 
associated 3-dimensional foliation)
extended this result to any
3-Sasakian dimension other than 11. Neither of the two authors noticed
Hitchin's result. Instead they considered geometry of the associated
foliation of $M$ by 3-dimensional leaves. It is likely that one could also
give a direct proof of Kashiwada's Theorem
working exclusively on $M.$ 
\bigskip
\centerline{\bf ACKNOWLEDGMENTS}
\medskip
Both authors would like to thank the Erwin Schr\"odinger
International Institute for Mathematical Physics in Vienna for
support and hospitality, where this work began. 
\bigskip\bigskip\bigskip

\bigskip
\bigskip
\medskip
\centerline{\bf Bibliography}
\medskip
\font\ninesl=cmsl9
\font\bsc=cmcsc10 at 10truept
\parskip=1.5truept
\baselineskip=11truept
\ninerm

\item{[BM1]} {\bsc A. Banyaga and P. Molino}, {\ninesl G\'eom\'etrie des formes
de contact compl\`etement integrable de type torique}, S\'eminaire Gaston
Darboux, Montpelier, 1991-1992, 1-25.
\item{[BM2]} {\bsc A. Banyaga and P. Molino}, {\ninesl Complete Integrability
in Contact Geometry}, Pennsylvania State University preprint PM 197, 1996.
\item{[Bes]} {\bsc A. L. Besse}, {\ninesl Einstein Manifolds}, 
Springer-Verlag, New York (1987). 
\item{[BG1]} {\bsc C. P. Boyer and  K. Galicki}, {\ninesl
3-Sasakian Manifolds}, 
to appear in {\it Essays on Einstein Manifolds}, International Press 1999;
C. LeBrun and M. Wang, Eds.
\item{[BG2]} {\bsc C. P. Boyer and  K. Galicki}, {\ninesl On Sasakian-Einstein
Geometry}, preprint DG/981108, October 1998.
\item{[BG3]} {\bsc C. P. Boyer and  K. Galicki}, {\ninesl A Note on Toric
Contact Geometry }, preprint DG/9907043 September 1999, to appear in The
Journal of Geometry and Physics.
\item{[Bl]} {\bsc D. E. Blair}, {\ninesl Contact Manifolds in
Riemannian Geometry}, Lecture Notes in Mathematics 509, Springer-Verlag,
New Yrok 1976.
\item{[BW]} {\bsc W.M. Boothby and H.C. Wang}, {\ninesl On Contact Manifolds},
Ann. of Math. 68 (1958), 721-734.
\item{[Car]} {\bsc Y. Carri\`ere}, {\ninesl Les propri\'et\'es topologiques
des flots riemanniens retrouv\'ees \`a l'aide du}
\break{{\ninesl th\'eor\`eme
des vari\'et\'es presque plates}, Math. Z. 186 (1984), 393-400.}
\item{[Gol]} {\bsc S.I. Goldberg}, {\ninesl Integrability of Almost Kaehler
Manifolds}, Proc. Amer. Math. Soc. 21 (1969), 96-100.
\item{[Gr]} {\bsc M. Gromov}, {\ninesl Volume and bounded cohomology}, Inst.
Hautes \'Etudes Sci. Publ. Math. 56 (1982), 213-307.
\item{[Hat]} {\bsc Y. Hatakeyama}, {\ninesl Some notes on differentiable
manifolds with almost contact structures}, T\^ohuku Math. J. 15 (1963),
176-181.
\item{[Hit]} {\bsc N.J. Hitchin}, {\ninesl The self-duality equations on a
Riemann surface}, Proc. London Math. Soc. 55 (1987), 59-126. 
\item{[IY]} {\bsc H. Inoue and K. Yano}, {\ninesl The Gromov invariant of
negatively curved manifolds}, Topology 21 (1981), 83-89.
\item{[Jel]} {\bsc W. Jelonek}, {\ninesl Positive and negative 
3-K-contact structures}, preprint, 1998.
\item{[Kas]} {\bsc T. Kashiwada}, {\ninesl On a contact 3-structure}, preprint.
\item{[LM]} {\bsc P. Libermann and C.-M. Marle}, {\ninesl Symplectic
Geometry and Analytical Mechanics}, D. Reidel Publishing Co., Dordrecht, 1987.
\item{[Mol1]} {\bsc P. Molino}, {\ninesl Riemannian Foliations}, Progress in
Mathematics 73, Birkh\"auser, Boston, 1988.
\item{[Mol2]} {\bsc P. Molino}, {\ninesl Feuilletages riemanniens sur les
vari\'et\'es compactes; champs de Killing transverses}, C.R. Ac. Sci. Paris
289 (1979), 421-423.
\item{[MoSe]} {\bsc P. Molino and Vlad Sergiescu}, {\ninesl Deux Remarques sur
Les Flots Riemanniens}, Manus. Math. 51 (1985), 145-161.
\item{[Rei1]} {\bsc B.L. Reinhart}, {\ninesl Foliated Manifolds with
Bundle-like Metrics}, Ann. of Math. 69 (1959), 119-132.
\item{[Rei2]} {\bsc B.L. Reinhart}, {\ninesl Differential Geometry of
Foliations}, Springer-Verlag, New York, 1983.
\item{[Ruk1]} {\bsc P. Rukimbira}, {\ninesl Chern-Hamilton's Conjecture and
K-Contactness}, Hous. J. of Math. 21 (1995), 709-718.
\item{[Ruk2]} {\bsc P. Rukimbira}, {\ninesl The Dimension of Leaf Closures of
K-Contact Flows}, Ann. Global Anal. and Geom. 12 (1994), 103-108.
\item{[Sek1]} {\bsc K. Sekigawa}, {\ninesl On some 4-dimensional compact
Einstein almost K\"ahler manifolds}, Math. Ann. 271 (1985), 333-337.
\item{[Sek2]} {\bsc K. Sekigawa}, {\ninesl On some compact
Einstein almost K\"ahler manifolds}, J. Math. Soc. Japan 39 (1987), 677-684.
\item{[Tan1]} {\bsc S. Tanno}, {\ninesl The automorphism groups of almost
contact Riemannian manifolds}, T\^ohoku Math. J. 21 (1969), 21-38.
\item{[Tan2]} {\bsc S. Tanno}, {\ninesl Geodesic flows
on $C_L$-manifolds and Einstein metrics on $S^3\times S^2$}, in
{\it Minimal submanifolds and geodesics (Proc.
Japan-United States Sem., Tokyo, 1977)}, pp. 283-292, North Holland,
Amsterdam-New York, 1979.
\item{[Tan3]} {\bsc S. Tanno}, {\ninesl Remarks on a 
Triple of K-contact Structures},
T\^ohoku Math. J. 48 (1996), 519-531.
\item{[Tho]} {\bsc C.B. Thomas}, {\ninesl Almost regular contact manifolds},
J. Differential Geom. 11 (1976), 521-533.
\item{[Ton]} {\bsc P. Tondeur}, {\ninesl Geometry of Foliations},
Birkh\"auser, Boston, 1997.
\item{[YK]} {\bsc K. Yano and M. Kon}, {\ninesl 
Structures on Manifolds}, Series in Pure Mathematics 3,  
World Scientific Pub. Co., Singapore, 1984.

\medskip
\bigskip \line{ Department of Mathematics and Statistics
\hfil November 1999} \line{ University of New Mexico \hfil} 
\line{ Albuquerque, NM 87131 \hfil } \line{ email: cboyer@math.unm.edu,
galicki@math.unm.edu\hfil} \line{ web pages:
http://www.math.unm.edu/$\tilde{\phantom{o}}$cboyer, 
http://www.math.unm.edu/$\tilde{\phantom{o}}$galicki \hfil}
\bye